\documentclass{article}
\usepackage[utf8]{inputenc}

\newtheorem{lemma}{Lemma}[section]

\newtheorem{theorem}[lemma]{Theorem}
\newtheorem{proposition}[lemma]{Proposition}
\newtheorem{corollary}[lemma]{Corollary}
\newtheorem{example}[lemma]{Example}

\newtheorem{note}[lemma]{Remark}

\title{Arities and aritizabilities of first-order theories\footnote{
The work of the author was carried out in the framework of the
State Contract of the Sobolev Institute of Mathematics, Project
No.~FWNF-2022-0012 (Section 1), and of Russian Scientific
Foundation, Project No.~22-21-00044 (Sections 2 and 3).}}
\author{Sergey V. Sudoplatov}
\date{}

\begin{document}

\maketitle
\begin{abstract}
We study and describe possibilities for arities of elementary
theories and of their expansions. Links for arities with respect
to Boolean algebras, to disjoint unions and to compositions of
structures are shown. Arities and aritizabilities are semantically
characterized. The dynamics for arities of theories is described.
\end{abstract}

{\bf Key words:} elementary theory, arity, expansion,
aritizability.

\bigskip
Arities of theories are important characteristics showing
complexity measures of theories \cite{ccegt96, SudTrans} and
reducing all definable sets to definable ones generated by
cylinders of special forms. It is closely linked with cylindric
algebras reflecting semantically first-order calculi
\cite{Maltsev, HMT1, HMT2, HMT3, ImLi}.

Special cases for arities of theories, especially binary, ternary
and related ordered theories are studied in a series of papers
including \cite{KulMac, Kulp1, Kulp2, Kulp3, Kulp4, Kulp5, Kulp6}.
Structures and links with respect to binary formulas are
investigated both in general case \cite{CCMCT, ShS, Su122, Erlag}
and for a series of natural classes of theories \cite{Su13cat,
Su13sm, KS15, EmIrk, EKS17, EmSud17, BEKPS, EKS19, EKS20}.

In the present paper we adapt the general cylindric approach and
describe semantically arities of theories, properties related to
the $n$-arity and $n$-aritizability of theories and their
dynamics.

The paper in organized as follows. In Section 1, we consider
arities of formulae and arities of theories, describe
possibilities of arities, describe arities for a series of natural
theories, characterize the $\omega$-categoricity and the stability
of $n$-ary theories. In Section 2, we introduce a series of
notions for the aritizability of a theory, describe sufficient
conditions and criteria for the aritizabilities in terms of
Boolean algebras, disjoint unions and compositions. In Section 3,
possibilities for aritizabilities are characterized semantically
and the dynamics for arities of theories is described.

Throughout we consider complete first-order theories $T$.

\section{$n$-ary formulae and theories}

{\bf Definition} \cite{ZPS}. A theory $T$ is said to be
\emph{$\Delta$-based}\index{Theory!$\Delta$-based}, where $\Delta$
is some set of~formulae without parameters, if any formula of $T$
is equivalent in $T$ to a~Boolean combination of formulae
in~$\Delta$.

For $\Delta$-based theories $T$, it is also said that $T$ has {\em
quantifier elimination}\index{Elimination of quantifiers} or {\em
quantifier reduction}\index{Reduction of quantifiers} up to
$\Delta$.

\medskip
{\bf Definition} \cite{CCMCT, ZPS}. {\rm Let $\Delta$ be a set of
formulae of a theory $T$, and $p(\bar{x})$ a type of $T$ lying in
$S(T)$. The type $p(\bar{x})$ is said to be
\emph{$\Delta$-based}\index{Type!$\Delta$-based} if $p(\bar{x})$
is isolated by a set of formulas $\varphi^\delta\in p$, where
$\varphi\in\Delta$, $\delta\in\{0,1\}$.}

\medskip
The following lemma, being a corollary of Compactness Theorem,
noticed in \cite{ZPS}.

\medskip
\begin{lemma}\label{lem01} A  theory  $T$  is
$\Delta$-based  if and only if,  for  any  tuple~$\bar{a}$ of any
{\rm (}some{\rm )} weakly saturated model of $T$, the type ${\rm
tp}(\bar{a})$ is $\Delta$-based.
\end{lemma}

\medskip
{\bf Definition} (cf. \cite{ccegt96}). An elementary theory $T$ is
called {\em unary}, or {\em $1$-ary}, if any $T$-formula
$\varphi(\overline{x})$ is $T$-equivalent to a Boolean combination
of $T$-formulas, each of which is of one free variable, and of
formulas of form $x\approx y$.

For a natural number $n\geq 1$, a formula $\varphi(\overline{x})$
of a theory $T$ is called {\em $n$-ary}, or an {\em $n$-formula},
if $\varphi(\overline{x})$ is $T$-equivalent to a Boolean
combination of $T$-formulas, each of which is of $n$ free
variables.

For a natural number $n\geq 2$, an elementary theory $T$ is called
{\em $n$-ary}, or an {\em $n$-theory}, if any $T$-formula
$\varphi(\overline{x})$ is $n$-ary.

A theory $T$ is called {\em binary} if $T$ is $2$-ary, it is
called {\em ternary} if $T$ is $3$-ary, etc.

We will admit the case $n=0$ for $n$-formulae
$\varphi(\overline{x})$. In such a case $\varphi(\overline{x})$ is
just $T$-equivalent to a sentence
$\forall\overline{x}\varphi(\overline{x})$.

If $T$ is a theory such that $T$ is $n$-ary and not $(n-1)$-ary
then the value $n$ is called the arity of $T$ and it is denoted by
${\rm ar}(T)$. If $T$ does not have any arity we put ${\rm
ar}(T)=\infty$.

Similarly, for a formula $\varphi$ of a theory $T$ we denote by
${\rm ar}_T(\varphi)$ the natural value $n$ if $\varphi$ is
$n$-ary and not $(n-1)$-ary. If $\varphi$ does not any arity we
put ${\rm ar}_T(\varphi)=\infty$. If a theory $T$ is fixed we
write ${\rm ar}(\varphi)$ instead of ${\rm ar}_T(\varphi)$.

\medskip
By the definition any $n$-theory is $\Delta_n$-based, where
$\Delta_n$ consists of formulae with $n$ free variables and
formulae of the form $x\approx y$. It implies that theories of
$n$-element models $\mathcal{M}$ are $n$-ary and based by formulae
describing these $n$-element structures and
differences/coincidences of elements.

\medskip
Using Lemma \ref{lem01} we obtain the following characterization
for the $n$-arity of a formula:

\begin{proposition}\label{pr02}
A $T$-formula $\varphi(\overline{x})$ is not $n$-ary if and only
if for any $T$-formulae $\psi_i(\overline{x}_i)$ with subtuples
$\overline{x}_i$ of the tuple $\overline{x}$ having
$l(\overline{x}_i)=n$ and
$T\vdash\varphi(\overline{x})\to\psi_i(\overline{x}_i)$, there
exists a tuple $\overline{a}\in\mathcal{M}\models T$ such that
$\mathcal{M}\models\psi_i(\overline{a}_i)\wedge\neg\varphi(\overline{a})$,
where $\overline{a}_i$ is a subtuple of $\overline{a}$ consisting
of substitutions of elements of $\overline{a}$ instead of
correspondent elements of $\overline{x}_i$.
\end{proposition}

By the definition the notion of $n$-arity is local and reduced to
finite sublanguages:

\begin{proposition}\label{pr03}
A theory $T$ of a language $\Sigma$ is $n$-ary if and only if for
any $T$-formula $\varphi(\overline{x})$ there is a finite
sublanguage $\Sigma'\subseteq\Sigma$ such that
$T\vdash\varphi(\overline{x})\leftrightarrow\psi(\overline{x})$,
where $\psi(\overline{x})$ is a Boolean combination of
$n$-formulae.
\end{proposition}

\begin{proposition}\label{pr031}
If $\mathcal{M}$ is a $n$-element structure, for $n\in\omega$,
then ${\rm ar}({\rm Th}(\mathcal{M}))\leq n$.
\end{proposition}

Proof. Since $|M|=n$ each ${\rm Th}(\mathcal{M})$-formula
$\varphi:=\varphi(x_1,\ldots,x_m)$ is ${\rm
Th}(\mathcal{M})$-equivalent to a disjunction of substitutions of
variables $x_{i_1},\ldots,x_{i_n}$ instead of $x_1,\ldots,x_m$
into the formula $\varphi$, as required.

\begin{note}\label{rem03} {\rm (cf. \cite{HMT1, ImLi})
Since negations of formulas with $n$ free variables again have $n$
free variables, witnessing the $n$-arity of a formula it suffices
to consider positive Boolean combinations of formulas with $n$
free variables, i.e., conjunctions and disjunctions of formulas
with $n$ free variables.

Thus for the description of definable sets for models
$\mathcal{M}$ of $n$-theories it suffices describe links between
definable sets $A$ and $B$ for $n$-formulas
$\varphi(\overline{x})$ and $\psi(\overline{y})$, respectively,
and definable sets $C$ and $D$ for
$\varphi(\overline{x})\wedge\psi(\overline{y})$ and
$\varphi(\overline{x})\vee\psi(\overline{y})$, respectively.

If $\overline{x}=\overline{y}$ then $C=A\cap B$ and $D=A\cup B$,
i.e., conjunctions and disjunctions work as set-theoretic
intersections and unions.

If $\overline{x}$ and $\overline{y}$ are disjoint then $C=A\times
B$ and
$D=(A+B)_{\mathcal{M}}\rightleftharpoons\{\langle\overline{a},\overline{b}\rangle\mid\overline{a}\in
A\mbox{ and }\overline{b}\in M, \mbox{ or } \overline{a}\in
M\mbox{ and }\overline{b}\in B\}$, i.e., $C$ is the Cartesian
product of $A$ and $B$, and $D$ is the ({\em generalized}) {\em
Cartesian sum} of $A$ and $B$ in the model $\mathcal{M}$.

If $\overline{x}\ne\overline{y}$, and $\overline{x}$ and
$\overline{y}$ have common variables, then $C$ and $D$ are
represented as a {\em mixed product} and a {\em mixed sum},
respectively, working partially as intersection and union, for
common variables, and partially as Cartesian product and Cartesian
sum, for disjoint variables.

If $\overline{x}$ and $\overline{y}$ consist of pairwise disjoint
variables and $\overline{x}\subseteq\overline{y}$ and
$\overline{x}\ne\overline{y}$ then for any formula
$\varphi(\overline{x})$ the set of solution of the formula
$\varphi(\overline{x})\wedge(\overline{y}\approx\overline{y})$ in
$\mathcal{M}$ is called a {\em cylinder} with respect to
$M^{l(\overline{y})}$ and generated by the set of solutions
$\varphi(\mathcal{M})$. In any case generating sets for cylinders
coincide their {\em projections}, i.e., sets of solutions for
formulas $\exists\overline{z}\varphi(\overline{x})$, where
$\overline{z}\subset\overline{x}$.

Since $n$-formulae produce cylinders on Cartesian products of
universes, definable sets of $n$-ary theories are composed by
Boolean combinations of definable cylinders, i.e., of elements of
cylindric algebras.}
\end{note}

{\bf Definition} (cf. \cite{SudTrans}). For a natural number $n$,
a theory $T$ is called {\rm $n$-transitive} if each $n$-type
$q(x_1,\ldots,x_n)\in S(T)$ is forced by its restriction to the
empty language.

\begin{proposition}\label{prTrans}
If a theory $T$ is $n$-transitive and non-$(n+1)$-transitive then
$T$ is not an $n$-theory.
\end{proposition}

Proof. Since $T$ is $n$-transitive, cylinders defined by
$n$-formulae are reduced to the cylinders defined by the formulae
for the empty language, i.e., they are defined by equalities and
inequalities. As $T$ is not $(n+1)$-transitive then there is a
$\emptyset$-definable set $X\subset M^{n+1}$ in a model
$\mathcal{M}\models T$ which in not reduced to the cylinders
defined by the formulae for the empty language. It means that a
formula $\varphi(x_1,\ldots,x_{n+1})$ defining $X$ is not
$T$-equivalent to $n$-formulae. Thus $T$ is not an $n$-theory, as
required.

\medskip
Clearly, generic constructions \cite{CCMCT, Su073} allow to
produce, for each $n\geq 1$, $n$-transitive and
non-$(n+1)$-transitive theories with unique $(n+1)$-ary predicates
and having quantifier elimination.

For instance, the theory $T$ of structure
$\mathcal{M}=\langle\{a,b,c,d\};R^{(3)}\rangle$ with the ternary
relation
$R=\{(a,b,c),(b,a,d),(b,c,d),(c,b,a),(a,c,d),(c,a,b),(c,d,a),$
$(d,c,b),(d,a,b),(a,d,c),$ $(b,d,a),(d,b,c)\}$ has quantifier
elimination, is $2$-transitive, not $3$-transitive, and thus ${\rm
ar}(T)=3$.

This example can be naturally spread for $n$-ary relations. In
view of Proposition \ref{prTrans} it implies the following:

\begin{corollary}\label{corTrans}
For any natural $n\geq 1$ there is a theory $T_n$ with ${\rm
ar}(T_n)=n$.
\end{corollary}

The following examples illustrate values ${\rm ar}(T)=n$.

\begin{example}\label{exunar}  {\rm \cite{SD1} For any theory $T_f$ of an unar, i.e., of one unary operation $f$,
${\rm ar}(T)\leq 2$. There are both theories $T_{f_1}$ with ${\rm
ar}(T_{f_1})=1$ and theories $T_{f_2}$ with ${\rm ar}(T_{f_2})=2$.
For instance, $f_1$ can be taken identical, and $f_2$~--- a
successor function on at least $3$-element set.}
\end{example}

\begin{example}\label{exacyclic}  {\rm \cite{SD1} For any theory $T_{\Gamma}$ of an acyclic graph $\Gamma$ with
unary predicates, ${\rm ar}(T_\Gamma)\leq 2$.}
\end{example}

\begin{example}\label{exsphere} {\rm
Let $E$ be the following equivalence relation on the set ${\bf
R}^{n}$:
$$\{(M,N)\mid M(x_1,\ldots,x_{n}), N(y_1,\ldots,y_{n})\in{\bf
R}^{n}, x^2_1+\ldots +x^2_{n}=y^2_1+\ldots +y^2_{n}\}.$$
Equivalence classes for the concentric spheres in ${\bf R}^{n}$
can not be reconstructed via cylinders defined by projections
which form concentric balls and circles. The homogeneity of
equivalence classes implies that each formula in the language
$\langle E \rangle$ is reduced to a Boolean combination of
$2n$-formulas. Thus ${\rm Th}(\langle{\bf R},E\rangle)$ is a
$2n$-theory which is not an $(2n-1)$-theory.

Adding a disjoint unary predicate $P$ and a bijection $f$ between
the set of spheres and $P$ we obtain names for spheres and an
additional coordinate for generating formulas for a basedness.
Thus we form a $(2n+1)$-theory which is not an $2n$-theory.

Hence all possibilities for ${\rm ar}(T)=n$ are realized.}
\end{example}

\begin{example}\label{exsurface} {\rm
Taking a non-degenerated algebraic surface at ${\bf R}^{n}$ which
is not reduced to cylinders we obtain a defining formula
$\varphi(\overline{x})$, $l(\overline{x})=n+1$, which is
$(n+1)$-formula and not an $n$-formula. In particular,
non-degenerated non-cylindrical surfaces of the second order in
${\bf R}^{3}$ are defined by formulas $\varphi$ with ${\rm
ar}(\varphi)=3$. For instance, taking the formula $x^2+y^2+z^2=1$
for the sphere $S$ we obtain projections $x^2+y^2\leq 1$,
$x^2+z^2\leq 1$, $y^2+z^2\leq 1$ which can not allow to
reconstruct $S$ by their Boolean combinations. }
\end{example}

\begin{example}\label{excircular} {\rm Recall \cite{KulMac, Kulp3, Kulp4}
that a {\em circular}, or {\em cyclic}  order relation is
described by a ternary relation $K_3$ satisfying the following conditions:\\
 (co1) $\forall x\forall y \forall z (K_3(x,y,z)\to K_3(y,z,x));$\\
 (co2) $\forall x\forall y \forall z (K_3(x,y,z)\land K_3(y,x,z)
 \leftrightarrow x=y \lor y=z \lor z=x);$\\
 (co3) $\forall x\forall y \forall z(K_3(x,y,z)\to \forall t[K_3(x,y,t)
 \lor K_3(t,y,z)]);$\\
 (co4) $\forall x\forall y \forall z (K_3(x,y,z)\lor K_3(y,x,z)).$

Clearly, ${\rm ar}(K_3(x,y,z))=3$ if the relation has at least
three element domain. Hence, theories with infinite circular order
relations are at least $3$-ary.

The following generalization of circular order produces a {\em
$n$-ball}, or {\em $n$-spherical}, or {\em $n$-circular} order
relation, for $n\geq 4$, which is described by a $n$-ary relation
$K_n$ satisfying the following conditions:\\
 (nbo1) $\forall x_1,\ldots,x_n (K_n(x_1,x_2,\ldots,x_n)\to K_n(x_2,\ldots,x_n,x_1));$\\
 (nbo2) $\forall x_1,\ldots,x_n \biggm(K_n(x_1,\ldots,x_i,x_{i+1},\ldots,x_n)\land$  $$\land K_n(x_1,\ldots,x_{i+1},x_{i},\ldots,x_n)
 \leftrightarrow\bigvee\limits_{i=1}^{n-1} x_i=x_{i+1}\biggm);$$\\
 (nbo3) $\forall x_1,\ldots,x_n(K_n(x_1,\ldots,x_n)\to \forall t[K_n(x_1,\ldots,x_{n-1},t)
 \lor K_n(t,x_2,\ldots,x_n)]);$\\
 (nbo4) $\forall x_1,\ldots,x_n (K_n(x_1,\ldots,x_i,x_{i+1},\ldots,x_n)\lor K_n(x_1,\ldots,x_{i+1},x_{i},\ldots,x_n)),$ $i<n.$

Clearly, ${\rm ar}(K_n(x_1,\ldots,x_n))=n$ if the relation has at
least $n$-element domain. Thus, theories with infinite $n$-ball
order relations are at least $n$-ary.}
\end{example}

\begin{theorem}\label{th04}
An $n$-ary theory $T$ is $\omega$-categorical if and only if there
are finitely many $T$-non-equivalent formulas with $n$ free
variables.
\end{theorem}

Proof. If $T$ is $\omega$-categorical then by Ryll-Nardzewski
Theorem there are finitely many $T$-non-equivalent formulas with
$m$ free variables for every $m$, in particular, for $m=n$.
Conversely, we again apply Ryll-Nardzewski Theorem showing that
are finitely many $T$-non-equivalent formulas with $m$ free
variables for every $m$. If $m\leq n$ then there are finitely many
$T$-non-equivalent formulas with $m$ free variables by the
monotony of this property with respect to the number of free
variables. If $m>n$ then by the $n$-arity of $T$ each $T$-formula
$\varphi(x_1,\ldots,x_m)$ is $T$-equivalent to a Boolean
combination of formulas with $n$ free variables. Since there are
finitely many $T$-non-equivalent possibilities for these formulas,
Boolean combinations produce finitely many possibilities, too. As
there are finitely many $T$-non-equivalent formulas with $m$ free
variables for every $m$ then $T$ is $\omega$-categorical by
Ryll-Nardzewski Theorem, as required.

\medskip
Recall \cite{Sh} that a formula $\varphi(\bar{x},\bar{y})$ of a
theory $T$ is {\em stable}\index{Formula!stable} if there are no
tuples $\bar{a}_n$, $\bar{b}_n$, $n\in\omega$, such that
$\models\varphi(\bar{a}_i,\bar{b}_j)\Leftrightarrow i\leq j$. The
theory $T$ is called {\em stable}\index{Theory!stable} if every
$T$-formula is stable.

In \cite{HaHa}, it was shown that any Boolean combination of
stable formulas is again a stable formula. Thus, using the
definition of $n$-ary theory we obtain the following:

\begin{theorem}\label{thstable}
An $n$-ary theory $T$ is stable if and only if each $T$-formula
with $n$ free variables is stable.
\end{theorem}

\section{Aritizable formulae and theories}

{\bf Definition.} A $T$-formula $\varphi(\overline{x})$ is called
{\em $n$-expansible}, or {\em $n$-arizable}, or {\em
$n$-aritizable}, if $T$ has an expansion $T'$ such that
$\varphi(\overline{x})$ is $T'$-equivalent to a Boolean
combination of $T'$-formulas with $n$ free variables.

A theory $T$ is called {\em $n$-expansible}, or {\em
$n$-arizable}, or {\em $n$-aritizable}, if there is an $n$-ary
expansion $T'$ of $T$.

A theory $T$ is called {\em arizable} or {\em aritizable}, if $T$
is $n$-aritizable for some $n$.

A $1$-aritizable theory is called {\em unary-able}, or {\em
unary-tizable}. A $2$-aritizable theory is called {\em
binary-tizable} or {\em binarizable}, a $3$-aritizable theory is
called {\em ternary-tizable} or {\em ternarizable}, etc.

\medskip
By the definition any $n$-theory is $n$-expansible, by itself, and
if $T$ is $n$-expansible then $T$ is $m$-expansible for each
$m>n$.

Besides each formula of an $n$-expansible theory is
$n$-expansible, too, but not vice versa in the following sense: if
each formula of a theory $T$ is $n$-expansible, it can not
guarantee that a resulting expansion $T'$, witnessing that
$n$-expansibility, is coordinated enough such that it is $n$-ary
or at least $n$-expansible.

\begin{proposition}\label{propbinfin} Any theory of a finite
structure $\mathcal{M}$ is binarizable.
\end{proposition}

Proof. Let $M=\{a_1,\ldots,a_m\}$. For any pair $\langle
a_i,a_j\rangle$, $i,j\leq m$, we introduce new binary singleton
predicate $B_{i,j}=\{\langle a_i,a_j\rangle\}$. We denote the
resulted expansion of $\mathcal{M}$ by $\mathcal{M}'$, and the
theory ${\rm Th}(\mathcal{M}')$ expanding given theory $T={\rm
Th}(\mathcal{M})$ by $T'$. Now an arbitrary $T'$-formula
$\varphi(\overline{x})$, with $\overline{x}=\langle
x_1,\ldots,x_n\rangle$, has finitely many solutions
$\overline{b}=\langle a_{k_1},\ldots,a_{k_n}\rangle$ in
$\mathcal{M}'$. We collect these solutions into a set $Z$. Without
loss of generality $Z\ne\emptyset$ since for  $Z=\emptyset$,
$T\vdash\varphi(\overline{x})\leftrightarrow\neg x_i\approx x_i$
for any $x_i$. Now the formula $\varphi(\overline{x})$ is
$T'$-equivalent to the following Boolean combination of binary
formulae:
$$
\bigvee_{\langle a_{k_1},\ldots,a_{k_n}\rangle\in
Z}\bigwedge_{i,j\leq n}B_{k_i,k_j}(x_i,x_j),
$$
as required.

\begin{note}\label{notefin} {\rm If $m,n\in\omega\setminus\{0\}$ and $M$ is an
$m$-element set, then $M^n$ has $2^{(m^n)}$ subsets producing
distinct $n$-ary predicates $Q_i$, $i<2^{(m^n)}$. Since by Stone
Theorem any finite Boolean has $2^l$ elements with $l$ generators,
there are $m^n$ independent predicates $Q'_j$ whose Boolean
combinations produce all these predicates. Taking a quantifier
free formula $\varphi(x_1,\ldots,x_k)$, for $k\geq n$, composed by
these independent predicates $Q'_j$ and having a perfect
disjunctive normal form we obtain $A^n_k\cdot 2^{m^n}$
possibilities for disjunctive members, where $A^n_k$ is used to
calculate the number of choice of $n$ variables among
$x_1,\ldots,x_k$ and there are $2^{m^n}$ possibilities for
positive and negative entries of $Q'_j$. Now there are
$2^{A^n_k\cdot 2^{m^n}}$ possibilities for
$\varphi(x_1,\ldots,x_k)$, big enough. At the same time, using the
arguments for Proposition \ref{propbinfin} we can obtain all
definable subsets of $M^n$ just using $m^2$ singleton binary
relations.}
\end{note}

Applying arguments for Proposition \ref{propbinfin} we immediately
obtain the following:

\begin{proposition}\label{propbinfinfor} Any formula of a theory having finitely many solutions is binarizable.
\end{proposition}

Proposition \ref{propbinfin} can be strengthened as follows:

\begin{proposition}\label{propunfin} Any theory of a finite
structure $\mathcal{M}$ is unary-tizable.
\end{proposition}

Proof. Let $M=\{a_1,\ldots,a_m\}$. For element $a_i$, $i\leq m$,
we introduce new unary singleton predicate $U_i=\{a_i\}$. We
denote the resulted expansion of $\mathcal{M}$ by $\mathcal{M}'$,
and the theory ${\rm Th}(\mathcal{M}')$ expanding given theory
$T={\rm Th}(\mathcal{M})$ by $T'$. Now an arbitrary $T'$-formula
$\varphi(\overline{x})$, with $\overline{x}=\langle
x_1,\ldots,x_n\rangle$, has finitely many solutions
$\overline{b}=\langle a_{k_1},\ldots,a_{k_n}\rangle$ in
$\mathcal{M}'$. We collect these solutions into a set $Z$. Without
loss of generality $Z\ne\emptyset$. Now the formula
$\varphi(\overline{x})$ is $T'$-equivalent to the following
Boolean combination of formulae each of which with one free
variable:
$$
\bigvee_{\langle a_{k_1},\ldots,a_{k_n}\rangle\in
Z}\bigwedge_{i,j\leq n}U_{k_i}(x_i),
$$
as required.

\medskip
Thus $|\mathcal{M}|$-many unary predicates produce a unary
expansion of the theory ${\rm Th}(\mathcal{M})$ with finite
$\mathcal{M}$. Besides using the proof of Proposition
\ref{propunfin} we have:

\begin{proposition}\label{propunfinfor} Any formula of a theory having finitely many solutions is unary-tizable.
\end{proposition}

\begin{note}\label{BA} {\rm By the definition for any natural $n$ both $n$-ary
formulae and $n$-aritizable formulae of a fixed theory $T$ are
closed under Boolean combinations. Therefore taking a model
$\mathcal{M}\models T$ and collecting in sets ${\rm
BA}_{kn}(\mathcal{M})$ and ${\rm BA}'_{kn}(\mathcal{M})$ definable
sets which are defined by $n$-ary, respectively, $n$-aritizable
formulae with $k$ free variables we obtain Boolean algebras
$\mathcal{BA}_{kn}(\mathcal{M})$ and
$\mathcal{BA}'_{kn}(\mathcal{M})$ of these definable sets.

Clearly,
$\mathcal{BA}_{kn}(\mathcal{M})\subseteq\mathcal{BA}'_{kn}(\mathcal{M})$,
and the equality
$\mathcal{BA}_{kn}(\mathcal{M})=\mathcal{BA}'_{kn}(\mathcal{M})$
means that any $n$-aritizable formula $\varphi(x_1,\ldots,x_n)$ of
$T$ is already $n$-ary.

In view of Proposition \ref{propunfinfor} the algebra ${\rm
BA}'_{kn}(\mathcal{M})$ satisfies the following condition: if $X$
and $Y$ are $\emptyset$-definable subsets of $M^k$ with finite
symmetric difference $X\div Y$ then $X\in{\rm
BA}'_{kn}(\mathcal{M})$ iff $Y\in{\rm BA}'_{kn}(\mathcal{M})$. At
the same time ${\rm BA}_{kn}(\mathcal{M})$ can be not closed under
finite symmetric difference since, for instance, there are
theories of finite structures which are not $n$-ary but by
Proposition \ref{propunfin} all theories of finite structures are
unary-tizable.

The Boolean algebras $\mathcal{BA}_{kn}(\mathcal{M})$ and
$\mathcal{BA}'_{kn}(\mathcal{M})$ have extensions
$\mathcal{BA}_{k}(\mathcal{M})$ and
$\mathcal{BA}'_{k}(\mathcal{M})$, respectively, consisting of
definable sets for $n$-ary/$n$-aritizable formulae with $m$ free
variables, for some $n$. Clearly, both
$\mathcal{BA}_{k}(\mathcal{M})$ and
$\mathcal{BA}'_{k}(\mathcal{M})$ equal the Boolean algebra
$\mathcal{B}_k(\mathcal{M})$ of all $\emptyset$-definable subsets
of $M^k$. Both these inclusions
$\mathcal{BA}_{kn}(\mathcal{M})\subseteq\mathcal{B}_k(\mathcal{M})$
and
$\mathcal{BA}'_{kn}(\mathcal{M})\subseteq\mathcal{B}_k(\mathcal{M})$
can be proper for $k>n$.

It was noticed above that aritizabilities of separated formulae of
a theory can not guarantee that witnesses of these aritizabilities
produce a $n$-ary theory, the possibility of {\em coordinated}
expansion of aritizable formulae is necessary. We denote by
$\mathcal{BA}''_{kn}(\mathcal{M})$ the Boolean algebra
$\mathcal{BA}'_{kn}(\mathcal{M})$ with a coordinated $n$-ary
expansion for all $n$-aritizable formulae with $k$ free
variables.}
\end{note}

Using Remark \ref{BA} we have the following characterizations of
$n$-arity and of $n$-aritizability of a theory $T$ in terms of
Boolean algebras of a model for $T$.

\begin{proposition}\label{BAar}
For any theory $T$, its model $\mathcal{M}$, and $n\in\omega$ the
following conditions hold:

$(1)$ $T$ is $n$-ary iff
$\mathcal{BA}_{kn}(\mathcal{M})=\mathcal{B}_k(\mathcal{M})$ for
each $k>n$;

$(1)$ $T$ is $n$-aritizable iff
$\mathcal{BA}''_{kn}(\mathcal{M})=\mathcal{B}_k(\mathcal{M})$ for
each $k>n$.
\end{proposition}

\medskip
{\bf Definition.} \cite{Wo} The {\em disjoint
union}\index{Disjoint union!of structures}
$\bigsqcup\limits_{n\in\omega}\mathcal{
M}_n$\index{$\bigsqcup\limits_{n\in\omega}\mathcal{ M}_n$} of
pairwise disjoint structures $\mathcal{ M}_n$ for pairwise
disjoint predicate languages $\Sigma_n$, $n\in\omega$, is the
structure of language
$\bigcup\limits_{n\in\omega}\Sigma_n\cup\{P^{(1)}_n\mid
n\in\omega\}$ with the universe $\bigsqcup\limits_{n\in\omega}
M_n$, $P_n=M_n$, and interpretations of predicate symbols in
$\Sigma_n$ coinciding with their interpretations in $\mathcal{
M}_n$, $n\in\omega$. The {\em disjoint union of
theories}\index{Disjoint union!of theories} $T_n$ for pairwise
disjoint languages $\Sigma_n$ accordingly, $n\in\omega$, is the
theory
$$\bigsqcup\limits_{n\in\omega}T_n\rightleftharpoons{\rm Th}\left(\bigsqcup\limits_{n\in\omega}\mathcal{ M}_n\right),$$
where\index{$\bigsqcup\limits_{n\in\omega}T_n$} $\mathcal{
M}_n\models T_n$, $n\in\omega$. Taking empty sets instead of some
structures $\mathcal{ M}_k$ we obtain disjoint unions of finitely
many structures and theories. In particular, we have the disjoint
unions $\mathcal{ M}_0\sqcup\ldots\sqcup\mathcal{ M}_n$ and their
theories $T_0\sqcup\ldots\sqcup T_n$.

\medskip
Clearly, disjoint unions of theories does not depend on choice of
correspondent disjoint unions of their models. Besides, disjoint
unions $\bigsqcup\limits_{n\in\omega}T_n$ are based by the unions
of the basing sets $\Delta_n$ for $T_n$ and by the formulae of the
form $P_n(x)$. Thus we have the following:

\begin{theorem}\label{thdu}
$1.$ For any theories $T_m$, $m\in\omega$, and their disjoint
union $\bigsqcup\limits_{m\in\omega}T_m$, all $T_m$ are
$n$-theories iff $\bigsqcup\limits_{m\in\omega}T_m$ is an
$n$-theory, moreover, ${\rm
ar}\left(\bigsqcup\limits_{m\in\omega}T_m\right)={\rm max}\{{\rm
ar}(T_m)\mid m\in\omega\}$.

$2.$ For any theories $T_m$, $m\in\omega$, and their disjoint
union $\bigsqcup\limits_{m\in\omega}T_m$, all $T_m$ are
$n$-aritizable iff $\bigsqcup\limits_{m\in\omega}T_m$ is
$n$-aritizable.
\end{theorem}

{\bf Definition} \cite{EKS20}. Let $\mathcal{M}$ and $\mathcal{N}$
be structures of relational languages $\Sigma_\mathcal{M}$ and
$\Sigma_\mathcal{N}$ respectively. We define the {\em composition}
$\mathcal{M}[\mathcal{N}]$ of $\mathcal{M}$ and $\mathcal{N}$
satisfying the following conditions:

1)
$\Sigma_{\mathcal{M}[\mathcal{N}]}=\Sigma_\mathcal{M}\cup\Sigma_\mathcal{N}$;

2) $M[N]=M\times N$, where $M[N]$, $M$, $N$ are universes of
$\mathcal{M}[\mathcal{N}]$, $\mathcal{M}$, and $\mathcal{N}$
respectively;

3) if $R\in\Sigma_\mathcal{M}\setminus\Sigma_\mathcal{N}$,
$\mu(R)=n$, then $((a_1,b_1),\ldots,(a_n,b_n))\in
R_{\mathcal{M}[\mathcal{N}]}$ if and only if $(a_1,\ldots,a_n)\in
R_{\mathcal{M}}$;

4) if $R\in\Sigma_\mathcal{N}\setminus\Sigma_\mathcal{M}$,
$\mu(R)=n$, then $((a_1,b_1),\ldots,(a_n,b_n))\in
R_{\mathcal{M}[\mathcal{N}]}$ if and only if $a_1=\ldots =a_n$ and
$(b_1,\ldots,b_n)\in R_{\mathcal{N}}$;

5) if $R\in\Sigma_\mathcal{M}\cap\Sigma_\mathcal{N}$, $\mu(R)=n$,
then $((a_1,b_1),\ldots,(a_n,b_n))\in
R_{\mathcal{M}[\mathcal{N}]}$ if and only if $(a_1,\ldots,a_n)\in
R_{\mathcal{M}}$, or $a_1=\ldots =a_n$ and $(b_1,\ldots,b_n)\in
R_{\mathcal{N}}$.

The composition $\mathcal{M}[\mathcal{N}]$ is called {\em
$e$-definable}, or {\em {\rm equ}-definable}, if
$\mathcal{M}[\mathcal{N}]$ has an $\emptyset$-definable
equivalence relation $E$ whose $E$-classes are universes of the
copies of $\mathcal{N}$ forming $\mathcal{M}[\mathcal{N}]$. If the
equivalence relation $E$ is fixed, the $e$-definable composition
is called {\em $E$-definable}.

\medskip
Using a nice basedness of $E$-definable compositions $T_1[T_2]$
(see \cite{EKS20}) till the formulas of form $E(x,y)$ and
generating formulas for $T_1$ and $T_2$ we have the following:

\begin{theorem}\label{thcomp} $1.$ For any theories $T_1$ and $T_2$ and
their $E$-definable composition $T_1[T_2]$, $T_1$ and $T_2$ are
$n$-theories, for $n\geq 2$, iff $T_1[T_2]$ is an $n$-theory,
moreover, ${\rm ar}(T_1[T_2])={\rm max}\{{\rm ar}(T_1),{\rm
ar}(T_2)\}$, if models of $T_1$ and of $T_2$ have at least two
elements, and ${\rm ar}(T_1[T_2])={\rm max}\{{\rm ar}(T_1),{\rm
ar}(T_2),2\}$, if a model of $T_1$ or $T_2$ is a singleton.

$2.$ For any theories $T_1$ and $T_2$ and their $E$-definable
composition $T_1[T_2]$, $T_1$ and $T_2$ are $n$-aritizable iff
$T_1[T_2]$ is $n$-aritizable.
\end{theorem}

Applying Proposition \ref{propunfin} and Theorem \ref{thcomp} we
immediately obtain:

\begin{corollary}\label{corcomp}
If each of theories $T_1$ and $T_2$ is a theory of a finite
structure, or of an infinite structure and $n$-arizable, then
their $E$-definable composition $T_1[T_2]$ is $n$-aritizable.
\end{corollary}

\section{Unary-tizable, binarizable and aritizable theories, their definable sets and dynamics}

Let $T$ be a theory with a unary expansion $T'$. Since unary
formulas $\varphi(x)$ and $\psi(y)$ have either equal or disjoint
free variables we do not have essential mixed sums and mixed
products forming definable sets for a model $\mathcal{M}$ of $T'$,
i.e., all definable sets are formed using unions, intersections,
Cartesian sums and Cartesian products of definable subsets of $M$,
without parameters.

Conversely, having a system of definable sets formed by unions,
intersections, Cartesian sums and Cartesian products of subsets of
$M$, we can introduce names for these subsets and generate, using
this introduced language, all given definable sets.

Thus we obtain the following characterization for the
unary-tizability of a theory in terms of definable sets:

\begin{theorem}\label{thunary}
A theory $T$ is unary-tizable if and only if for any {\rm
(}some{\rm )} model $\mathcal{M}$ of $T$ any definable set is
formed by unions, intersections, Cartesian sums and Cartesian
products of subsets of $M$.
\end{theorem}

Similarly, all definable sets of binarizable theories are
generated by unions, intersections, Cartesian sums and Cartesian
products of subsets of $M^2$, extended by mixed sums and mixed
products of these subsets and their combinations:

\begin{theorem}\label{thbinary} A theory $T$ is binarizable if and only if
for any {\rm (}some{\rm )} model $\mathcal{M}$ of $T$ any
$\emptyset$-definable set is formed by unions, intersections,
Cartesian sums, Cartesian products, mixed sums and mixed products
of subsets of $M^2$.
\end{theorem}

By Theorem \ref{thbinary} definable sets of binarizable theories
are generated by combinations of 3-dimensional cylinders with
two-dimensional generators.

Theorems \ref{thunary} and \ref{thbinary} admit the following
natural generalizations based on $(n+1)$-dimensional cylinders
with $n$-dimensional generators.

\begin{theorem}\label{th_n_ary} A theory $T$ is $n$-aritizable, for $n\geq
1$, if and only if for any {\rm (}some{\rm )} model $\mathcal{M}$
of $T$ any $\emptyset$-definable set is formed by unions,
intersections, Cartesian sums, Cartesian products, mixed sums and
mixed products of subsets of $M^n$.
\end{theorem}

\begin{theorem}\label{th_ary}
A theory $T$ is aritizable if and only if for any {\rm (}some{\rm
)} model $\mathcal{M}$ of $T$ any $\emptyset$-definable set is
formed by unions, intersections, Cartesian sums, Cartesian
products, mixed sums and mixed products of subsets of $M^n$, for
some $n$.
\end{theorem}

\begin{note}\label{remar} {\rm Using Theorem \ref{th_n_ary} one can form
a definable subset $X\subset M^k$ of an infinite model
$\mathcal{M}$ of a $n$-theory, for $k>n$, such that $X$ has an
infinite complement and each projection of $X$ equals $M^m$ for
some $m<n$. It implies that Boolean combinations of these
projections can not reconstruct $X$. Thus the structure $\langle
M, X\rangle$ has a $n$-expansible $k$-theory.}
\end{note}

{\bf Definition.} A $T$-formula $\varphi(\overline{x})$ is called
{\em constantizable} if $T$ has an expansion $T'$ such that
$\varphi(\overline{x})$ is $T'$-equivalent to a Boolean
combination of formulae of forms $x\approx y$ and $x\approx c$
with variables $x$, $y$ and constants $c$, i.e., $T$ has an
expansion $T''$ such that $\varphi(\overline{x})$ is
$T'$-equivalent to a Boolean combination of formulae of forms
$x\approx y$ and of formulae of unary singleton predicates whose
solutions consist of constants.

A theory $T$ is called {\em constantizable} if any $T$-formula
$\varphi(\overline{x})$ is constantizable.

\medskip
By the definition any constantizable theory is unary-tizable, but
not vice versa, as the following assertions show.

\begin{proposition}\label{propconstfor} A formula $\varphi(\overline{x})$ of a theory $T$ is constantizable iff $\varphi(\overline{x})$ is $T$-equivalent to a Boolean
combination of formulae of form $x\approx y$ and formulae with
finitely many solutions.
\end{proposition}

Proof. Let $\varphi(\overline{x})$ be a constantizable formula.
Since $\varphi(\overline{x})$ is unary-tizable we can divide
$\varphi(\overline{x})$ in some expansion of $T$ onto cases with
distinct/equal values for free variables. Thus without loss of
generality $\varphi(\overline{x})$ is a Boolean combination of
formulae $(x_i\approx x_j)^\delta$ and $\psi(x_i)$, for
$x_i,x_j\in\overline{x}$, $i\ne j$, $\delta\in\{0,1\}$, written in
a disjunctive normal form. Since $\varphi(\overline{x})$ is
constantizable, its definable set $A$ in a model of $T$ is
represented by a Boolean combination of cylinders for $x\approx y$
and $x\approx c$. Thus, $\varphi(\overline{x})$ is $T$-equivalent
to a Boolean combination of formulae of form $x\approx y$ and
formulae with finitely many solutions.

Now let $\varphi(\overline{x})$ is $T$-equivalent to a Boolean
combination of formulae of form $x\approx y$ and formulae
$\psi(\overline{x})$ with finitely many solutions. We may assume
that $\varphi(\overline{x})$ consistent and represented as a
disjunctive normal form. We collect in a set $Z$ all finite sets
of solutions for formulae $\psi(\overline{x})$. Now for the finite
set $\cup Z$ we apply the construction for Proposition
\ref{propunfin} reducing the set $Z'$ of all coordinates for
tuples in $\cup Z$ to unary singleton predicates and so to the
formulae $x\approx c$, for $c\in Z'$. It implies that a Boolean
combination of these formulae and formulae of the form $x\approx
y$ is equivalent to $\varphi(\overline{x})$, i.e.,
$\varphi(\overline{x})$ is constantizable, as required.

\medskip
Proposition \ref{propconstfor} immediately implies:

\begin{corollary}\label{corconst1} A theory $T$ is constantizable iff each $T$-formula is $T$-equivalent to a Boolean
combination of formulae of form $x\approx y$ and formulae with
finitely many solutions.
\end{corollary}

\medskip
{\bf Definition} \cite{BaLa}. A theory $T$ is called {\em strongly
minimal}\index{Theory!strongly minimal} if for any formula
$\varphi(x,\bar{a})$ of language obtained by adding parameters of
$\bar{a}$ (in some model ${\cal M}\models T$) to the language of
$T$, either $\varphi(x,\bar{a})$, or $\neg\varphi(x,\bar{a})$ has
finitely many solutions.

\medskip
Using Corollary \ref{corconst1} we obtain:

\begin{corollary}\label{corconst2} Any constantizable theory is strongly minimal.
\end{corollary}

Now we consider some dynamics of arities of theories under
expansions. Since the property of non-$n$-arizability forbids
$n$-ary expansions it suffices to study possibilities for
expansions of $n$-aritizable theories.

\begin{proposition}\label{propnarim}
A theory $T$ has a non-aritizable expansion iff $T$ has an
infinite model.
\end{proposition}

Proof. If $T$ has an infinite model there are expansions $T'$ of
$T$ collecting, for instance, examples \ref{exsphere},
\ref{exsurface}, \ref{excircular} forbidding $n$-arity for each
$n$. Thus $T'$ is not aritizable.

Conversely, if $T$ has an finite model then each expansion $T'$ of
$T$ has a finite model producing aritizability of $T'$ by
Proposition \ref{propunfin}, as required.

\medskip
Using examples above we observe that for each natural $n\geq 1$
there are theories $T_{n}$ with ${\rm ar}(T_{n})=n$ and finite
models. Thus there are theories $T_{kn}$ with ${\rm ar}(T_{kn})=k$
and ${\rm ar}(T'_{kn})=n$ for some expansions $T'_{kn}$ of
$T_{kn}$.

Besides, for each natural $n\geq 1$ there are:

1) theories $T_{n,\infty}$ with ${\rm ar}(T_{n,\infty})=n$ and
${\rm ar}(T'_{n,\infty})=\infty$ for some expansions
$T'_{n,\infty}$ of $T_{n,\infty}$: it suffices to expand a $n$-ary
theory with infinite models by new predicates forbidding the
$k$-aritizability for each $k>n$;

2) theories $T_{\infty,n}$ with ${\rm ar}(T_{\infty,n})=\infty$
and ${\rm ar}(T'_{\infty,n})=n$ for some expansions
$T'_{\infty,n}$ of $T_{\infty,n}$: it suffices to expand an
$n$-aritizable theory which is not $m$-ary for any $m$ till a
$n$-theory.

Thus the arities can be freely increased and decreased and we
obtain the following:

\begin{theorem}\label{tharim}
For any $\mu,\nu\in(\omega\setminus\{0\})\cup\{\infty\}$ there is
a theory $T_{\mu,\nu}$ and its expansion $T'_{\mu,\nu}$ such that
${\rm ar}(T_{\mu,\nu})=\mu$ and ${\rm ar}(T'_{\mu,\nu})=\nu$.
\end{theorem}

\section{Conclusion}

We considered possibilities for arities of theories and their
dynamics, reductions of formulas to ones of special forms as well
as definable sets connected with these reductions. It can be used
both for databases, simplifying them to ones with bounded
dimensions, for geometric objects represented as finite
combinations of cylinders, and for cryptographic constructions
representing complicated configurations by simpler ones. It would
be interesting to describe values of arities and aritizabilities
for natural classes of theories.

\noindent Sobolev Institute of Mathematics, \\ 4, Acad. Koptyug
avenue, Novosibirsk, 630090, Russia; \\ Novosibirsk State
Technical University, \\ 20, K.Marx avenue, Novosibirsk, 630073,
Russia

\medskip\noindent
e-mail: sudoplat@math.nsc.ru

\end{document}